\documentclass[12pt]{amsart} %{sn-jnl}
\usepackage{amsmath,amssymb,amsfonts,amsthm,amstext,verbatim,url}
\usepackage{graphicx}

\RequirePackage[colorlinks,citecolor=blue,urlcolor=blue]{hyperref}
\usepackage{breakurl}

\theoremstyle{plain}
\newtheorem{theorem}{Theorem}

\newtheorem{conjecture}[theorem]{Conjecture}
\newtheorem{question}[theorem]{Question}

\numberwithin{equation}{section}
\numberwithin{figure}{section}
\numberwithin{theorem}{section}

\theoremstyle{remark}

\newcounter{mycount}
\newenvironment{romlist}{\begin{list}{\rm(\roman{mycount})}%
   {\usecounter{mycount}\labelwidth=1cm\itemsep 0pt}}{\end{list}}

\newenvironment{letlist}{\begin{list}{\rm(\alph{mycount})}%
   {\usecounter{mycount}\labelwidth=1cm\itemsep 0pt}}{\end{list}}

\newcommand\eps{\epsilon}

\newcommand\oo{\infty}

\newcommand\HH{{\mathbb H}}

\newcommand\QQ{{\mathbb Q}}

\newcommand\sQ{{\mathcal Q}}

\newcommand\sF{{\mathcal F}}
\newcommand\ZZ{{\mathbb Z}}

\newcommand\RR{{\mathbb R}}
\newcommand\PP{{\mathbb P}}

\newcommand\wt{\widetilde}

\newcommand\sm{\setminus}
\renewcommand\a{\alpha}

\newcommand\Si{\Sigma}
\newcommand\si{\sigma}

\newcommand\g{\gamma}
\newcommand\resp{respectively}

\renewcommand\o{\text{\rm o}}

\newcommand\qq{\qquad}

\newcommand\SLE{\text{\rm SLE}}
\newcommand\ec{\eps_{\text{\rm c}}}
\newcommand\lra{\leftrightarrow}

\newcommand\vpc{\vec p_{\text{\rm c}}}
\newcommand\what{\widehat}
\newcommand\ac{\alpha_{\text{\rm c}}}
\newcommand\var{\text{\rm var}}
\newcommand\nlra{\nleftrightarrow}
\newcommand\covid{{\textsc{Covid}}-19}

\title{Selected problems in probability theory}

\author{Geoffrey R.\ Grimmett}
\address{Statistical Laboratory, Centre for
Mathematical Sciences, Cambridge University, Wilberforce Road,
Cambridge CB3 0WB, UK} 
\email{grg@statslab.cam.ac.uk}
\urladdr{\url{http://www.statslab.cam.ac.uk/~grg/}}

\date{9 April 2022} %\today}
\dedicatory{Dedicated in friendship  to Catriona Byrne}
\begin{document}
\begin{abstract}
This celebratory article contains a personal and idiosyncratic selection of a few open problems
in discrete probability theory.   These include 
certain well known questions concerning Lorentz scatterers and self-avoiding walks, 
and also some problems of percolation-type. The author hopes the reader will find something to
leaven winter evenings, and perhaps even a project for the longer term. 
  
\end{abstract}


\maketitle

\section*{Personal remarks}

The editorial team of Springer Mathematics has become almost family for many of us worldwide,
with Catriona Byrne at its heart. 
She has come to know us better than we know ourselves, always with sympathy, and
with an honest and constructive approach to occasionally challenging areas of professional debate. 
Through our numerous collaborations, she and I have kindled a warm friendship 
that will persist into the next phase of our adventures. We wish her many happy years 
free from the woes of authors, editors, readers, and publishers.

\section{Introduction}\label{sec:intro}

Probability has been a source of many tantalising problems over the centuries,
of which the St Petersburg paradox, Fermat's problem of the points, and Bertrand's random triangles
feature still in introductory courses.  Whereas the conceptual problems of the past are now largely resolved,
contemporary questions arise frequently where the intuitive apparatus of sub-fields collide.
Many prominent problems are to be found at the conjunction of probability and 
discrete geometry.
This short and idiosyncratic article summarises some of these. This account is personal and incomplete, 
and is to be viewed as a complete review of nothing. 
The bibliography is not intended to be complete, and apologies are
extended
to those whose work has been omitted.

The questions highlighted here vary from the intriguing to the profound. Whereas some may seem like puzzles with
limited consequence, others will require new machinery and may have far-reaching implications.   

The problem of counting self-avoiding walks
is introduced in Section \ref{sec:saw},
with emphasis on the existence of critical exponents and the scaling limit in two 
dimensions, followed by a fundamental counting problem on a random percolation cluster.
Section \ref{sec:Lor} is devoted to Lorentz scatterers and the Ehrenfest wind/tree model, 
followed by Poissonian mirrors in two dimensions, and finally Manhattan pinball.
Two well-known conjectures concerning product measures 
are presented in Section \ref{sec:ineq},
namely the bunkbed conjecture and the negative association of a uniform spanning forest.
Section \ref{sec:rosl} is concerned with the identification of criticality for the
randomly oriented square lattice. In the final Section \ref{sec:dse}, we present
two basic problems associated with a model for a dynamic spatial epidemic,
provoked in part by \covid.

\section{Self-avoiding walks}\label{sec:saw}

\subsection{Origins}
Self-avoiding walks were first introduced in the chemical theory of 
polymerisation (see \cite{f,Orr}),
and their properties have received much attention since from mathematicians and physicists 
(see, for example, \cite{bdgs,GZ-rev,ms}). 

A path in an infinite graph $G=(V,E)$ is called \emph{self-avoiding} if no vertex is visited more than once. 
Fix a vertex $v\in V$, and let $\Si_n(v)$ be the set of $n$-step self-avoiding walks (SAWs) starting at $v$.
The principal combinatorial problem is to determine how the cardinality $\si_n(v):=|\Si_n(v)|$ 
grows as $n\to\oo$, and the complementary
probability problem is to establish properties of the shape of a  randomly selected member
of $\Si_n(v)$. Progress has been striking but limited.

\subsection{Asymptotics}

It is now regarded as elementary that the so-called \emph{connective constant} $\kappa=\kappa(G)$, given by
$$
\log\kappa = \lim_{n\to\oo} \frac 1n \log\si_n(v),
$$
exists when $G$ is quasi-transitive, and is
independent of the choice of $v$. Thus, in this case
$$
\si_n(v) = \kappa^{(1+\o(1))n}.
$$
The correction term is much harder to understand. We shall not make
precise the concept of a $d$-dimensional lattice,
but for definiteness the reader may concentrate  on the hypercubic lattice $\ZZ^d$.

\begin{conjecture}\label{conj-saw}
For $d\ge 2$ there exists a \emph{critical exponent} $\g=\g_d$ such that the following holds.
Let $G$ be a $d$-dimensional lattice. There exists a constant $A>0$ such that\footnote{A 
logarithmic correction is in fact expected in \eqref{eq:corr} when $d=4$.}
\begin{equation}\label{eq:corr}
\si_n \sim A n^{\g-1}\kappa(G)^n \qq\text{as } n\to\oo.
\end{equation}
Furthermore,
$$
\g=\begin{cases} \frac{43}{32} &\text{when } d=2,\\
1 &\text{when } d \ge 4.
\end{cases}
$$
\end{conjecture}

 See \cite{bdgs,ms} for further discussion and results so far, and the papers  \cite{HS1,HS2}
 of Hara and Slade when $G=\ZZ^d$ with $d\ge 5$, for which case they prove that $\g=1$.
Of particular interest is the case when $G=\HH$, the hexagonal lattice.  
By a beautiful exact calculation that verifies an earlier conjecture of Nienhuis \cite{Nien}
based in conformal field theory, Duminil-Copin and Smirnov \cite{ds} proved that
$$
\kappa(\HH) = \sqrt{1+ \sqrt 2}.
$$
The proof reveals a discrete holomorphic function that is highly suggestive of a connection to 
a Schramm--Loewner evolution (see \cite{Kemp}), namely the following. 

\begin{question}
Does a uniformly distributed  $n$-step SAW from the origin of $\HH$ converge weakly,
when suitably rescaled,
to the Schramm--Loewner random curve $\SLE_{8/3}$?
\end{question}

Progress on this question should come hand-in-hand with a calculation of the
associated critical exponent $\gamma=\frac{43}{32}$. Gwynne and Miller \cite{GM21} have proved
the corresponding weak limit in the universe of Liouville quantum gravity.

\subsection{Self-avoiding walks in a random environment}
How does the sequence $(\si_n)$ behave when the underlying graph $G$ is random?
For concreteness,
we consider here the infinite cluster $I$ of bond percolation on $\ZZ^2$ with edge-density $p>\frac12$
(see \cite{G99}). 

\begin{question}
Does the limit $\mu(v):= \lim_{n\to\oo} \si_n(v)^{1/n}$
exist a.s., and satisfy $\mu(v)=\mu(w)$ a.s.\ on the event $\{v,w \in I\}$?
\end{question}

Related discussion, including of the issue of deciding when $\mu(v)=p\mu(\ZZ^2)$ a.s.\ 
on the event $\{v \in I\}$, may be found in papers of Lacoin \cite{Lac1, Lac2}.
The easier SAW problem on (deterministic) weighted graphs is considered in \cite{GL-wtsaw}. 

\section{Lorentz scatterers}\label{sec:Lor}

\subsection{Background}
The scattering problem of Lorentz \cite{Lor} gives rise to the following general question.
Scatterers are distributed randomly about $\RR^d$. Light is shone from the origin in a given direction,
and is subjected to reflection at the scatterers.
Under what circumstances is the light ray: (i) bounded, (ii) unbounded, (iii) diffusive?  
While certain special cases are understood, the general question remains open.
The problems mentioned here are concerned with \emph{aperiodic} 
distributions of scatterers;
the periodic case is rather easier.

\begin{figure}
\centerline{\includegraphics[width=0.5\textwidth]{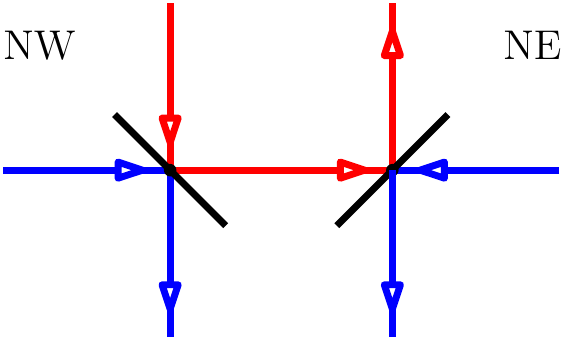}}
   \caption{A NW and a NE mirror. Each is reflective on both sides.}
\label{fig:mirrors}
\end{figure}

\subsection{Ehrenfest wind/tree model}\label{ssec3.2}
The following notorious problem on the square lattice $\ZZ^2$ has resisted solution for many years.
Let $p\in[0,1]$.
At each vertex of $\ZZ^2$ is placed a mirror with probability $p$, or alternatively nothing.
Mirrors are plane and two-sided. 
Each mirror is designated a \emph{north-east} (NE) mirror with probability $\frac12$,
or alternatively a \emph{north-west} (NW) mirror. The states of different vertices are independent.
The meanings of the mirrors are illustrated in Figure \ref{fig:mirrors}.

Light is shone from the origin in a given compass direction, 
say north, and it is reflected off the surface of any mirror encountered.
The problem is to decide whether or not the light ray is unbounded.

\begin{question}\label{q3.1}
Let $\theta(p)$ be the probability that the light ray is unbounded. 
For what values of $p$ is it the case that $\theta(p)>0$?
\end{question}

\begin{figure}
\centerline{\includegraphics[width=0.6\textwidth]{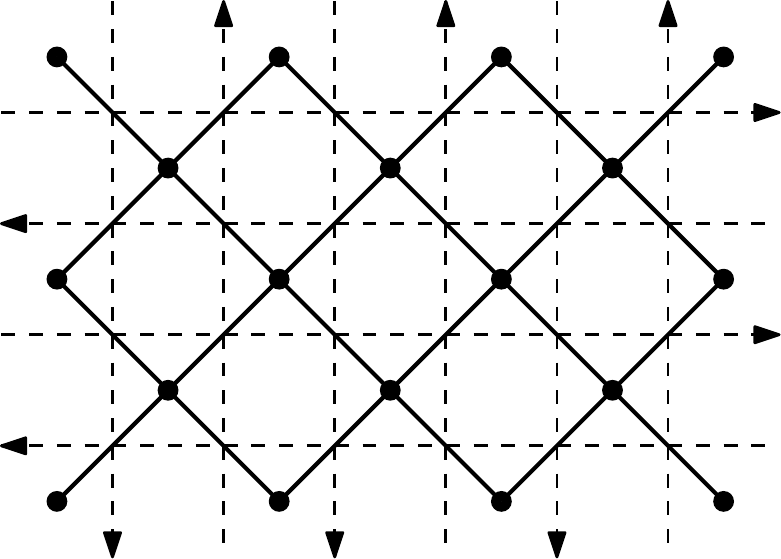}}
   \caption{From the original (dashed) square lattice $\ZZ^2$ one may construct a diagonal lattice $\what\ZZ^2$.
In fact there are two such diagonal lattices, and this fact may be used to 
obtain some information about the power-law behaviour of the light ray when $p=1$.
The Manhattan orientations are not relevant to the usual Ehrenfest model, but
are provided to facilitate the discussion of Manhattan pinball in Section \ref{sec:mp}.}
\label{fig:lordiag}
\end{figure}

It is trivial that $\theta(0)=1$. By considering bond percolation 
(with density $p/2$) on the diagonal lattice of
Figure \ref{fig:lordiag}, and using the fact that there is no percolation
when $p=1$, one obtains the less trivial fact that $\theta(1)=0$ (see \cite{G99}).
Very little more is known rigorously about the answer to Question \ref{q3.1}.

\subsection{Poisson mirrors}
\begin{figure}
\centerline{\includegraphics[width=0.5\textwidth]{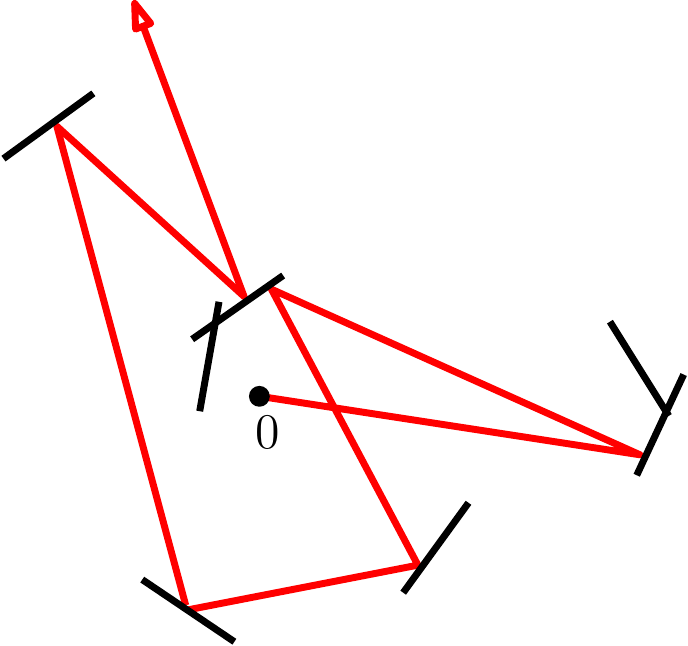}}
   \caption{Light from the origin is reflected off the mirrors.}
\label{fig:poiss}
\end{figure}

Here is version of the wind/tree model in the two-dimensional continuum $\RR^2$. Let $\Pi$ be a rate-$1$ Poisson
process in $\RR^2$. Let $\eps>0$, and let $\mu$ be a probability measure on $[0,\pi)$. We possess
an infinity of two-sided, plane mirrors of length $\eps$, and we centre one at each point in $\Pi$; the inclination 
to the horizontal of each mirror is random with law $\mu$, and different mirrors have independent inclinations.
Think of a mirror as being a randomly positioned, closed line segment of length $\eps$, and
let $M$ denote the union of these segments.
We call $\mu$ \emph{degenerate} if it is concentrated on a single atom, and shall assume
$\mu$ is non-degenerate. See Figure \ref{fig:poiss}.

Light is shone from the origin at an angle $\alpha$ to the horizontal. Let $I_\a$ be the indicator function
of the event that the light ray is unbounded.
Some convention is adopted for the zero-probability event that the light strikes  an intersection of two or more mirrors. 

We may assume that the origin $0$ does not lie in $M$.
Let $C$ be the component of $\RR^2\sm M$ containing $0$, and let $\{0\lra\oo\}$ be the event that
$C$ is unbounded.
It is a standard result of so-called needle percolation (see \cite{Roy}) that
there exists $\ec=\ec(\mu)\in(0,\oo)$ such that
$$
\PP_\mu(0\lra\oo) \begin{cases} >0 &\text{when } \eps<\ec,\\
=0 &\text{when } \eps>\ec.
\end{cases}
$$
(Here and later, the subscript $\mu$ keeps track of the choice of $\mu$.) 
Obviously, on the event that $0\nlra\oo$, we have that $I_\alpha=0$ for all $\alpha$.
Therefore, 
\begin{equation*}
\PP_\mu(I_\alpha=1\text{ for some } \alpha)=0, \qq \eps>\ec.
\end{equation*}
The converse issue is much harder and largely open.

\begin{question}\label{qn3.2}
Suppose $\mu$ is non-degenerate. 
\begin{letlist}
\item
Does there exist $\ec'=\ec'(\mu)>0$ such that $\theta_\mu(\eps)>0$ for $\eps<\ec'$?
\item
Could it be that $\ec'(\mu)=\ec(\mu)$?
\item
In particular, what happens when $\mu$ is the uniform measure on $[0,\pi)$?
\end{letlist} 
\end{question}

Let $\sQ$ be the set of probability measures $\mu$ that are non-degenerate and have support in 
the rational angles $\pi\QQ$. Suppose $\mu\in\sQ$ and $0 <\eps < \ec(\mu)$.
Harris \cite{Har} has shown the striking fact that, $\PP_\mu$-a.s.\ on the event that $0\lra\oo$, we have
that $I_\alpha=1$ for (Lebesgue) almost every $\alpha$. 

This leads to a deterministic question. Let $K$ be the set of mirror configurations for which
$0\lra\oo$ but $I_\alpha=0$ for all $\alpha$.

\begin{question}
Is $K$ non-empty?
\end{question}

Harris's theorem implies in effect that  $\PP_\mu(K)=0$ when $\mu\in\sQ$ and $\eps\ne \ec(\mu)$.
Question \ref{qn3.2}(c) hints at the possibility that
 $\PP_\mu(K)=0$ when $\mu$ is the uniform measure on $[0,\pi)$ and $\eps<\eps(\mu)$.

Here is  a final question concerning diffusivity.
Let $\mu\in\sQ$, and denote by
$X_\alpha(t)$ the position at time $t$
of the light ray that leaves the origin at angle $\alpha$.

\begin{question}
Is it the case that, on the event $I_\alpha$, $X_\alpha(\cdot)$ is diffusive?
That is, the limit $\si^2:=\lim_{t\to\oo} t^{-1}\var(X_\alpha(t))$ exists in $(0,\oo)$,
and, when normalized, $X_\alpha(t)$ is asymptotically normally distributed.
\end{question}

Related work on Lorentz models in the so-called Boltzmann--Grad limit may be found in \cite{LT1,LT2}.

\subsection{Manhattan pinball}\label{sec:mp}
Here is a variant of the Ehrenfest model motivated by a problem of quantum localization,
\cite[Sec.\ 4.2]{Car}, \cite[p.\ 238]{Spen}.
Draw $\ZZ^2$ and the diagonal lattice $\what\ZZ^2$ as in Figure \ref{fig:lordiag};
each edge of $\ZZ^2$ receives its Manhattan orientation as indicated in the figure.
Consider bond percolation with density $q$ on the diagonal lattice.
Along each open edge of $\what\ZZ^2$ we place a two-sided plane mirror.  Light is shone from the origin along
a given one of the two admissible directions, and it is reflected by any mirror that it encounters
(such reflections are automatically consistent with the Manhattan orientations).
Let $\theta(q)$ be the probability that the light ray is unbounded. 

\begin{question}
Could it be that $\theta(q)=0$ for all $q>0$?
\end{question}

It follows as in Section \ref{ssec3.2} that $\theta(q)=0$ for $q\ge \frac12$,
and it has been proved by Li in \cite{li} that there exists $\eps>0$
such that $\theta(q)=0$ when $q>\frac12-\eps$. The proof uses the method of enhancements; 
see \cite{AG}, \cite[Sect.\ 3.3]{G99}.

\section{Two stochastic inequalities}\label{sec:ineq}

\subsection{Bunkbed inequality}

The mysterious \lq bunkbed' inequality was posed by Kasteleyn around 1985 (see \cite[Rem.\ 5]{vdBK01}, and also \cite{OH98}).
Of its various flavours, we  select the following. Let $G=(V,E)$ be a finite simple graph.
From $G$ we construct two copies denoted $G_1=(V_1,E_1)$ and $G_2=(V_2,E_2)$.
For $v\in V$ we write $v_i$ for the copy of $v$ lying in $V_i$.
We now attach $G_1$ and $G_2$ 
by adding edges $\langle v_1,v_2\rangle$ for each $v\in V$. 
This new graph is denoted $\wt G$, and it may be
considered as the product graph $G \times K_2$ where $K_2$ is the complete graph on two vertices (that is, an edge).
We may think of the $G_i$ as being \lq horizontal' and the extra edges as being \lq vertical'.

Each edge of $\wt G$ is declared \emph{open} with probability $p$,
independently of the states of other edges. Write $\PP_p$ for the appropriate product measure. For two vertices $u_i,v_j$ of
$\wt G$, we write $\{u_i \lra v_j\}$ for the event that there exists a $u_i/v_j$ path using
only open edges.

\begin{conjecture}
For $u,v\in V$, we have 
$$
\PP_p(u_1\lra v_1) \ge \PP_p(u_1\lra v_2).
$$
\end{conjecture}

There is uncertainty over whether this was the exact conjecture of Kasteleyn.
For example, it is suggested in \cite{HKN} (and perhaps elsewhere) that Kasteleyn
may have made the stronger conjecture  that
the inequality holds even after conditioning on the set $T$ of open vertical edges.
 
Some special cases of the bunkbed conjecture have been proved (see the references in \cite{HKN},
and more recently \cite{TR}), 
but the general question remains open.

 \subsection{Negative correlation}

Our next problem is quite longstanding (see \cite{Pem}) and remains mysterious.
In a nutshell it is to prove that the uniform random forest measure (USF) 
has a property of  negative association.

Let $G=(V,E)$ be a finite graph which, for simplicity, we assume has neither loops nor multiple edges.
A subset $F\subseteq E$ is called a \emph{forest} if $(V,F)$ has no cycles. Let $\sF$ be the
set of all forests in $G$ and let $\Phi$ be a random forest chosen uniformly from $\sF$. We call $\Phi$ 
\emph{edge-negatively associated} if
\begin{equation}\label{eq:ena}
\PP(e,f\in \Phi) \le \PP(e\in \Phi) \PP(f \in \Phi), \qq e,f\in E,\ e\ne f.
\end{equation}

\begin{conjecture}\label{conj:usf}
For all graphs $G$, the random forest $\Phi$ is edge-negatively associated.
\end{conjecture}

One may formulate various forms of negative dependence, amongst which 
the edge-negative association of \eqref{eq:ena} is quite weak.
One may 
conjecture that $\Phi$ has a stronger variety of such dependence. 
Further discussion may be found  in \cite[Sec.\ 3.9]{G-rcm}
and \cite{Pem}.

Experimental evidence for Conjecture \ref{conj:usf} is quite strong. A similar conjecture may be made for 
uniform measure on the set of $F\subseteq E$ such that $(V,F)$ is \emph{connected} (abbreviated to UCS). 
In contrast, uniform spanning tree (UST) is
well understood via the Kirchhoff theory of electrical networks, and further by \cite{FM}.
USF, UCS, and UST are special cases of the so-called random-cluster measure with cluster-weighting factor
$q$ satisfying $q<1$ (see \cite[Sects 1.5, 3.9]{G-rcm}). 

In recent work, \cite{BCH, BCHS}, the percolative properties of the 
weighted random forest (or \lq arboreal gas')
on $\ZZ^d$ have been explored. It turns out
that there is a phase transition if and only if $d\ge 3$. 

\section{Randomly oriented square lattice}\label{sec:rosl}

The following percolation-type problem remains open.
Consider the square lattice $\ZZ^2$ and let $p\in[0,1]$. 
 Each horizontal edge is oriented rightward with probability $p$, and otherwise leftward. Each vertical edge is oriented
 upward with probability $p$, and otherwise downward.  Write $\vec\ZZ^2$ for the ensuing randomly oriented network.
 
 Let $\theta(p)$ denote the 
 probability that the origin $0$ is the endpoint of an infinite path of $\vec\ZZ^2$ 
 that is oriented away from $0$. The challenge is to determine for which $p$ it is 
 the case that $\theta(p)>0$. It is elementary that $\theta(0)=1$, and that
 $\theta(p)=\theta(1-p)$. It is less obvious that $\theta(\frac12)=0$ (see \cite[1st edn]{G99}), which is proved via a coupling with bond percolation. By a comparison 
 with oriented percolation, we have that $\theta(p)>0$ if $p>\vpc$, 
 where $\vpc$ is the critical point of oriented percolation on $\ZZ^2$; it is not difficult
 to deduce by the enhancement method (see \cite{AG}, \cite[Sect.\ 3.3]{G99})
 that there exists $p'\in(\frac12,\vpc)$ such 
 that $\theta(p)>0$ when $p>p'$. It is believed that $\vpc \sim 0.64$, and proved that $\vpc<0.6735$.

\begin{question}
Is it the case that $\theta(p)>0$ for $p\ne\frac12$.
\end{question}

It is shown in \cite{GG01} that, for all $p$, $\vec\ZZ^2$
 is either critical or supercritical in the following sense: if any small
positive density of oriented edges is added at random, then there is a strictly
positive probability that the origin is the endpoint of an infinite self-avoiding
oriented path in the resulting graph.

\section{Dynamic stochastic epidemics}\label{sec:dse}

The recent pandemic has inspired a number of mathematical problems, including
the following stochastic model (see \cite{GZ22}). Particles are placed at time $0$ at the points
of a rate-$1$ Poisson process in $\RR^d$, where $d \ge 1$. Each particle diffuses around 
$\RR^d$ according to a Brownian motion, 
independently of other particles. At any given time,
each particle is in one of the states S (susceptible), I (infected), R (removed/dead).

At time $0$ there exists a unique particle in state I, and all others are in state S.
The infection/removal rules are as follows. 
\begin{letlist}
\item If an infected particle comes within distance
$1$ of a susceptible particle, the latter particle is infected. 
\item An infected particle remains infected for a period of time having the exponential distribution with parameter $\alpha$,
and is then removed. 
\end{letlist}
We call this the \lq diffusion model'.

\emph{Survival} is said to occur if, with a strictly positive probability, 
infinitely many particle are ultimately infected. It is proved in
\cite{GZ22} that, when $d \ge 1$ and $\alpha$ is sufficiently large, survival does not occur.
The following two questions (amongst others) are left open.

\begin{question}\label{qn:3}\mbox{\hfil}
\begin{romlist}
\item When $d=1$, could it be that there is no survival for any $\alpha>0$?
\item When $d\ge 2$, does survival occur for sufficiently small $\alpha>0$?
\end{romlist}
\end{question}

In a variant of this problem termed the \lq delayed diffusion model', 
a much fuller picture is known. Suppose, instead of the above, 
a particle moves only when it is infected; susceptible particles are stationary.
The answers to Question \ref{qn:3}(i,\,ii) are then no and yes, \resp, and indeed (when $d\ge 2$) there
exists a critical value $\ac(d)\in(0,\oo)$ of $\alpha$ marking the onset of survival. The key difference between
the two systems is that the latter model has a property of monotonicity that is lacking 
in the former.

The delayed diffusion model is a continuous space/time cousin of the discrete-time \lq frog model' of 
\cite{AMP02a} (see also \cite{RS02}),
with the addition with removal. Further relevant references may be found in \cite{GZ22}.

\providecommand{\bysame}{\leavevmode\hbox to3em{\hrulefill}\thinspace}
\providecommand{\MR}{\relax\ifhmode\unskip\space\fi MR }
% \MRhref is called by the amsart/book/proc definition of \MR.
\providecommand{\MRhref}[2]{%
  \href{http://www.ams.org/mathscinet-getitem?mr=#1}{#2}
}
\providecommand{\href}[2]{#2}

%\bibliography{cb1}
%\bibliographystyle{amsplain}

\end{document}